\documentclass{article}

\usepackage{amsmath,amssymb}
\usepackage{amsthm}
\usepackage{palatino}
\usepackage{epsfig}
\usepackage{graphicx}

\newcommand{\la}{\lambda}
\newcommand{\be}{\beta}
\newcommand{\al}{\alpha}
\newcommand{\ga}{\gamma}
\newcommand{\ep}{\epsilon}
\newcommand{\et}{\eta}
\newcommand{\de}{\delta}
\newcommand{\ka}{\kappa}
\newcommand{\mcl}{\mathcal{L}}
\newcommand{\lone}{|\mcl_{\Roman{counter}}|}
\newcommand{\Lra}{\Leftrightarrow}

\author{Claire Levaillant}

\title{A quantum combinatorial approach for computing a tetrahedral network of Jones-Wenzl projectors}
\date{January 5, 2013}
\begin{document}
\maketitle
\begin{center}
Abstract
\end{center}
Trivalent plane graphs are used in various areas of mathematics
which relate for instance to the colored Jones polynomial, invariants
of 3-manifolds and quantum computation. Their evaluation is based on computations in the Temperley-Lieb algebra and more specifically the Jones-Wenzl projectors. We use the work in \cite{KL} to
present a quantum combinatorial approach for evaluating a tetrahedral net. 
On
the way we recover two equivalent definitions for the unsigned Stirling
numbers of the first kind and we provide an equality for the quantized
factorial using these numbers.

\section{Introduction and Background}
This paper presents a combinatorial heuristic approach for computing in the generic case a tetrahedral network of Jones-Wenzl projectors, inspired by \cite{KL} but different from \cite{KL}. Tetrahedral networks have a wide range of applications which include quantum computation and invariants of $3$-manifolds. In \cite{BL}, they are used to construct unitary braid groups representations. In \cite{GA}, they are used for finding a formula for the colored Jones function of the simplest hyperbolic non-$2$-bridge knot. The Jones-Wenzl projectors on top of their own algebraic interest constitute a mathematical model associated with the trivalent vertex which represents the fusion of two particules in anyonic systems, see \cite{WA}. They first appear in \cite{WE}. A both recent and broad reference on Jones-Wenzl projectors and their categorification is \cite{CS}.
 Our paper is structured as follows. We begin by setting some definitions and notations and by recalling some elementary facts. We then proceed to a heuristic evaluation of the tetrahedron, followed by a discussion on a generalization of our heuristics to any trivalent plane graph, thus also enlightening the difference between our approach and the one of \cite{KL}. In the last section, we provide a formula for the quantized factorial which uses unsigned Stirling numbers of the first kind.
\subsection{Smoothings of a crossing and the bracket polynomial}

In what follows, $d$ is a complex number and $A^2$ and $A^{-2}$ are the two complex roots of the quadratics $X^2+dX+1$, so that $d=-A^2-A^{-2}$. \\
For a crossing, we speak of an $A$-smoothing (resp $A^{-1}$-smoothing) when the two regions swept out by turning the over-crossing line counter-clockwise is labeled by $A$ (resp $A^{-1}$), like on the figure below. \begin{center}
\epsfig{file=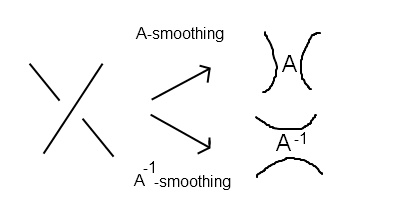, height=4cm}
\end{center}
We recall below the definition of the Kauffman bracket polynomial for $K$ an unoriented link.
$$<K>\,=\sum_{all\; possible\; smoothings} product\;of\;the\;labels\times d^{\;number\; of\; Jordan\; curves}$$
With this definition, we obviously have the following identity\\
\epsfig{file=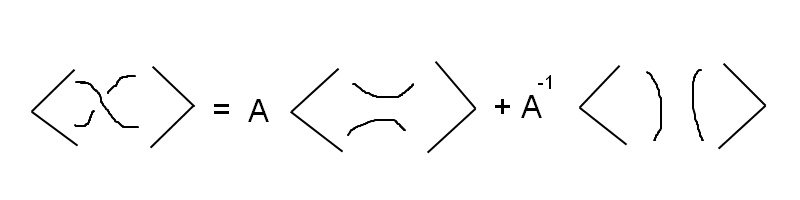, height=3cm}\\
where each small diagram must be seen as part of a larger one.
\subsection{Temperley-Lieb algebras, Jones-Wenzl projectors and their recursion formulas}
In what follows, $T_n$ will denote the $n$-strand Temperley-Lieb algebra over $\mathbb{C}$, with generators $1_n$, $U_1$, $\dots$, $U_{n-1}$ and relations
\begin{eqnarray*}
U_i^2&=&d\,U_i\\
U_iU_{i\pm 1}U_i&=&U_i\\
U_iU_j&=&U_jU_i\qquad if\qquad |i-j|>1
\end{eqnarray*}
In terms of diagrams, $U_i$ is the $n$-tangle that has two horizontal strands joining nodes $i$ and $i+1$, one at the top and one at the bottom, and no intersection between the vertical strands. Multiplying two $U_i$'s is made by stacking one diagram on top of the other and following the strands. The complex variable $d$ is called the loop value.
For instance $U_iU_{i+1}$ is the following $n$-tangle.
\begin{center}
\epsfig{file=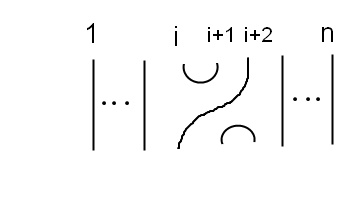, height=3cm}
\end{center}
The Jones-Wenzl projectors $f_0$, $f_1$, $\dots$, $f_{n-1}$ are $n$ elements of the Temperley-Lieb algebra $T_n$ that are recursively defined by
\begin{eqnarray*}
f_0&=&1_n\\
f_{i+1}&=&f_i-\mu_{i+1}\;\;f_i\,U_{i+1}\,f_i
\end{eqnarray*}
with
$$\mu_1=\frac{1}{d}\;\;\text{and}\;\;
\mu_{i+1}=\frac{1}{d-\mu_i}$$

It is a straightforward verification that is left to the reader that
\begin{eqnarray*}
f_i^2&=&f_i\;\;\;\;\;\;\;\;\;\;\;\;\;\;\;\forall\;i\\
f_i\,U_j&=&U_j\,f_i\,=\,0\;\;\;\forall j\leq i
\end{eqnarray*}

Conversely, we have the following lemma.

\newtheorem{Lemma}{Lemma}
\begin{Lemma}
If $g$ is a non-zero element of $T_n$ which satisfies
\begin{eqnarray*}
g^2&=&g\\
g\,U_i&=&U_i\,g\,=\,0\;\;\;\forall 1\leq i\leq n-1,
\end{eqnarray*}
then $g=f_{n-1}$
\end{Lemma}

\noindent The proof of this lemma appears in \cite{KL} and we briefly recall it here. If $g$ is only a linear combination of terms containing turn-backs, then by the first equality of the lemma, we can write $g$ as a product of $g$ times a linear combination of terms containing turn-backs. By the second equality of the lemma, such a product is zero, hence $g$ would be zero, a contradiction. Then $g$ must be a linear combination of $1_n$ and of the $U_i$'s with the coefficient of $1_n$ being non-zero. Since $g^2=g$, this coefficient must in fact be $1$. So, $$g=1_n +\mathcal{U},$$ where $\mathcal{U}$ is a linear combination of the $U_i$'s.\\
Notice the projector $f_{n-1}$ satisfies all the hypotheses of the lemma, hence we may also write $$f_{n-1}=1_n+\mathcal{U^{'}}$$ By writing the product $g\,f_{n-1}$ in two different ways, namely
\begin{equation*}\begin{split}g\,f_{n-1}&=g(1_n+\mathcal{U^{'}})=g\\
&=(1_n+\mathcal{U})\,f_{n-1}=f_{n-1},\end{split}\end{equation*}
one obtains $g=f_{n-1}$.\hfill $\square$ \\
We use the lemma to give an explicit definition of the Jones-Wenzl projectors, one that does not use recursion. Like in \cite{KL}, we define an element of $T_n$ as follows.
\newtheorem{Definition}{Definition}
\begin{Definition}
$$\underset{|}{\overset{\!\!\!\!n|}{\boxed{\begin{array}{l}\\\end{array}}}}=\frac{1}{\lbrace n\rbrace !}\,\sum_{\sigma\in S_n} (A^{-3})^{t(\sigma)}\;\bigg<\underset{|}{\overset{|}{\boxed{\begin{array}{l}\hat{\sigma}\end{array}}}}\bigg>$$
with $$\begin{array}{l}
\ast\;t(\sigma)\;\text{is the number of transpositions in the minimal decomposition of $\sigma$}\\ \text{as a product of transpositions $(i,i+1)$.}\\
\ast\;\hat{\sigma}\;\text{is the $n$-tangle obtained by replacing each transposition $(i,i+1)$}\\\text{in the minimal decomposition of $\sigma$ by a braid.}\\
\ast\;\begin{array}{l}\lbrace n\rbrace !=\sum_{\sigma\in S_n} (A^{-4})^{t(\sigma)}=\prod_{k=1}^n\frac{1-(A^{-4})^k}{1-A^{-4}}\end{array}
\end{array}$$
\end{Definition}
\noindent On the right hand side of the equality, it is to be understood that we expand the tangle inside the bracket by using the bracket identity from the end of $\S\,1.1$. By doing so, we obtain an element in the Temperley-Lieb algebra. For the reader interested in having explicit formulas for the projectors, we would like to mention that Scott Morrison in \cite{SM} gives a neat method for calculating the coefficients in front of the Temperley-Lieb diagrams in the expansion of the Jones-Wenzl projectors. It is shown in \cite{KL} that the element of Definition $1$ satisfies the conditions of Lemma $1$, thus this element is precisely the Jones-Wenzl projector $f_{n-1}$.
\begin{Definition}
We define $\Delta_n$ as the bracket evaluation of the closure of $f_{n-1}$ (\emph{i.e} the trace of $f_{n-1}$), which we simply write
\begin{center}
\epsfig{file=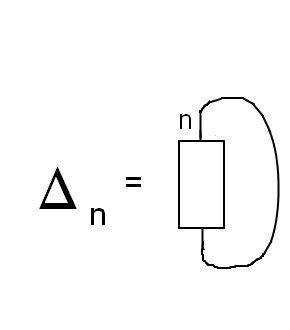, height=4cm}
\end{center}
\end{Definition}
\noindent We now use the recursive defining formula for the Jones-Wenzl projector $f_n$ in order to find a recursive formula satisfied by $\Delta_n$. We have, where we first consider the closure only on the last strand.
\begin{center}\epsfig{file=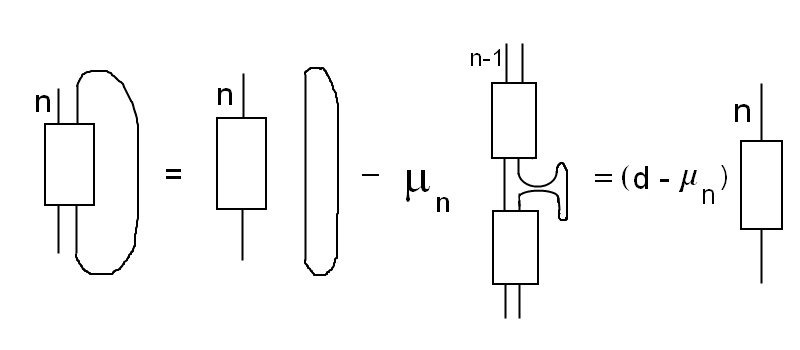, height=5cm}\end{center}
\noindent By now taking the closure on the first $n$ strands and doing a bracket evaluation, we get
\begin{equation}\Delta_{n+1}=(d-\mu_n)\,\Delta_n\end{equation} In particular, we derive
\begin{center}
\epsfig{file=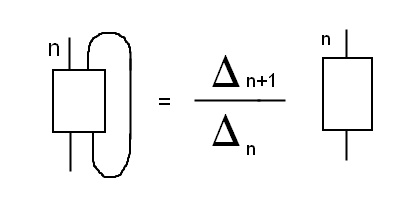, height=4cm}
\end{center}
and also \begin{equation}\mu_{n+1}=\frac{\Delta_n}{\Delta_{n+1}}\end{equation}
Gathering Eq. $(1)$ and Eq. $(2)$, we obtain the recursion
$$\left\lbrace\begin{array}{l}
\Delta_{n+2}=d\,\Delta_{n+1}-\Delta_n\\
\Delta_0=1\;\;\Delta_1=d
\end{array}\right.$$
This is a recursion formula for the Chebyshev polynomial of the second kind. The recursion is easily solved by replacing $d$ with $-A^2-A^{-2}$, which yields
$$\Delta_{n+2}+A^{-2}\,\Delta_{n+1}=-A^2(\Delta_{n+1}+A^{-2}\,\Delta_n).$$
We now derive
$$\Delta_n=(-1)^n\,\frac{A^{2n+2}-A^{-2n-2}}{A^2-A^{-2}}$$
Letting $q=A^2$, we get \begin{equation}\Delta_n=(-1)^n\,[n+1]_q,\end{equation} where $[n]_q$ is the quantum integer defined by
$$[n]_q=\frac{q^n-q^{-n}}{q-q^{-1}}$$
After recalling these basic facts about the Jones-Wenzl projectors, we define the trivalent vertex in terms of Jones-Wenzl projectors. By definition,
\begin{center}\epsfig{file=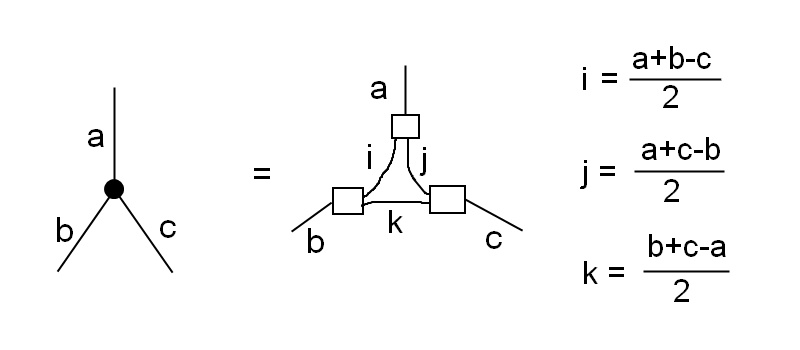, height=4cm}\end{center}
\noindent We impose that $a+b-c$, $a+c-b$ and $b+c-a$ are all positive and even. Our goal in the next section is to propose a method for evaluating the tetrahedron
\begin{center}
\epsfig{file=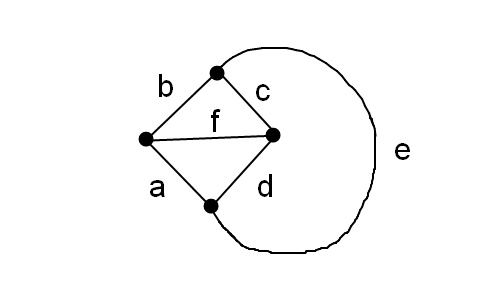, height=4cm}
\end{center}
\section{Evaluation of the tetrahedron}
\subsection{Preliminary detour}
Notice $A=-1$ is convenient. Indeed, since by definition,
$$\lbrace n\rbrace= 1+A^{-4}+\dots + (A^{-4})^{n-1},$$
we see that $$\lbrace n\rbrace (A=-1)=n$$
Further, when $A=-1$, the quantity $(A^{-3})^{t(\sigma)}$ is simply the signature $\varepsilon(\sigma)$ of the permutation $\sigma$. So, $f_{n-1}$ simplifies to
\begin{center}
\epsfig{file=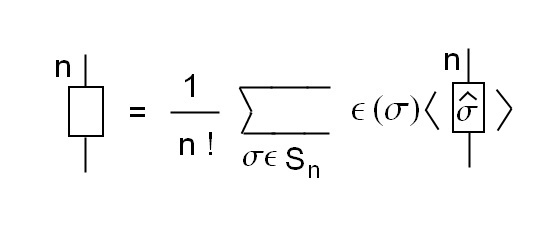, height=4cm}
\end{center}
Furthermore, when $A=-1$, we have the nice properties that loops can be removed by the bracket polynomial and that braids are simply elements of the symmetric group.
\begin{center}
\epsfig{file=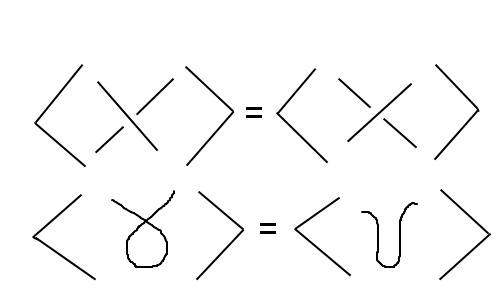, height=4cm}
\end{center}
Let $V$ be a vector space over $\mathbb{C}$ with $dim\,V=N$ and let $(e_1,\dots,\,e_{N})$ be a basis of $V$. Given $n\leq N$ and $\sigma\in\,S_n$, an idea of \cite{KL} is to define a map
$$\hat{\sigma}:\,\begin{array}{ccc}
V^{\otimes n}&\longrightarrow& V^{\otimes n}\\
e_{i_1}\otimes\dots\otimes e_{i_n}&\longmapsto& e_{i_{\sigma(1)}}\otimes\dots\otimes e_{i_{\sigma(n)}}
\end{array}$$
and by analogy with the formula above for the Jones-Wenzl projector $f_{n-1}$ taken at $A=-1$ to let
$$P_n^N=\frac{1}{n!}\sum_{\sigma\in S_n}\varepsilon(\sigma)\,\hat{\sigma}$$
The map $\hat{\sigma}$ permutes the basis vectors in the tensor product, following the permutation $\sigma$. Given integers $i_1$, $i_2$, $\dots$, $i_n$ and $j_1$, $j_2$, $\dots$, $j_n$, chosen among $N$ integers and not necessarily distinct, let $i=(i_1,\,i_2,\,\dots,\,i_n)$ and $j=(j_1,\,j_2,\,\dots,\,j_n)$. Then, the coefficient $(i,j)$ in the matrix representing the linear transformation $\hat{\sigma}$ is given by
$$[\hat{\sigma}]_{j_1\dots j_n}^{i_1\dots i_n}=\delta_{j_{\sigma(1)}}^{i_1}\,\delta_{j_{\sigma(2)}}^{i_2}\dots\delta_{j_{\sigma(n)}}^{i_n}$$
Consider the Jones-Wenzl projector defined like above when $A$ is specialized to the value $-1$, which we will denote by $f_{n-1}(A=-1)$. Then, we have the unexpected and quite nice following result.
\newtheorem{Theorem}{Theorem}
\begin{Theorem}
The trace of the Jones-Wenzl projector $f_{n-1}(A=-1)$ evaluated at $d=N$ is the trace of the linear transformation $P_n^{N}$. In other words, when $A=-1$, we have $\Delta_n(d=N)=Tr(P_n^N)$.
\end{Theorem}
\noindent\textbf{Proof.} The first trace is defined by closing the $n$-tangle inside the bracket polynomial with loops that match the outgoing lines with the respective incoming lines at the same position. The second trace is the usual standard trace. \\
On the one hand, visibly the trace of $f_{n-1}(A=-1)$ evaluated at $d=N$ is
$$\frac{1}{n!}\,\sum_{\sigma\in S_n} \varepsilon(\sigma)\,N^{s(\sigma)},$$
where $s(\sigma)$ denotes the number of cycles in the decomposition of $\sigma$ into a product of cycles with disjoint support. \\
On the other hand, the trace of $P_n^N$ is given by
$$\frac{1}{n!}\,\sum_{\sigma\in S_n}\,\varepsilon(\sigma)\,\sum_{a_1\dots a_n}[\hat{\sigma}]_{a_1\dots a_n}^{a_1\dots a_n}$$
The coefficient inside the second sum is zero, unless
$$\left\lbrace\begin{array}{ccc}
a_{\sigma(1)}&=&a_1\\
a_{\sigma(2)}&=&a_2\\
&\vdots&\\
a_{\sigma(n)}&=&a_n
\end{array}\right.$$
This means that the $a_k$'s are constant on the support of each cycle composing $\sigma$. Like above, if $s(\sigma)$ denotes the number of cycles of the permutation $\sigma$, there are $N^{s(\sigma)}$ possible choices for the $a_i$'s. Then,
\begin{equation}Tr(P_n^N)=\frac{1}{n!}\,\sum_{\sigma\in S_n}\,\varepsilon(\sigma)\,N^{s(\sigma)}\end{equation}
So we see that both traces are equal, as claimed in the theorem.\hfill$\square$\\
We further give a second way of evaluating the trace of the linear transformation $P_n^N$. By computing this trace in two different ways, we will derive a nice combinatorial equality. The key remark is the following.
\begin{Lemma}
If $a_i=a_j$ with $i\neq j$, then $[P_n^N]_{b_1\dots b_n}^{a_1\dots a_n}=0$.
\end{Lemma}
\noindent \textbf{Proof.} Suppose $a_i=a_j$ with $i\neq j$. For future reference, denote by $\tau$ the transposition $(i,j)$. We have $$[P_n^N]_{b_1\dots b_n}^{a_1\dots a_n}=\frac{1}{n!}\,\sum_{\sigma\in\mathcal{I}}\,\varepsilon(\sigma)$$
where $\mathcal{I}=\lbrace \sigma\in S_n,\;\forall k\in\lbrace 1,\dots,n\rbrace,\; b_{\sigma(k)}=a_k\rbrace$. We have
$$\left|\begin{array}{l}
b_{\sigma\tau(i)}=b_{\sigma(j)}=a_j=a_i\\
b_{\sigma\tau(j)}=b_{\sigma(i)}=a_i=a_j\\
b_{\sigma\tau(k)}=b_{\sigma(k)}=a_k\;\;\;\;\;\;\;\;\;\text{when $k\not\in\lbrace i,j\rbrace$}
\end{array}\right.$$
Thus, the map
$$\begin{array}{ccc} \mathcal{I}&\longrightarrow & \mathcal{I}\\
\sigma&\longmapsto&\sigma\circ\tau\end{array}$$
is well defined. Moreover, this map is an involution. It follows that
$$\sum_{\sigma\in\mathcal{I}}\,\varepsilon(\sigma)=\sum_{\sigma\in\mathcal{I}}\,\varepsilon(\sigma\circ\tau)=0,$$
where the last equality holds since for $\tau$ transposition, $\varepsilon(\sigma\circ\tau)=-\varepsilon(\sigma)$ so that the second sum is the opposite of the first one. \hfill$\square$\\
By the Lemma, we now have
$$Tr(P_n^N)=\sum_{\begin{array}{l} a_1,\dots,a_n\\\text{all distinct}\end{array}}\,[P_n^N]_{a_1\dots a_n}^{a_1\dots a_n}$$
Further, when the $a_k$'s are all distinct, the value of $[\hat{\sigma}]_{a_1\dots a_n}^{a_1\dots a_n}$ is one only when $\sigma$ is the identity and the signature of the identity is one. There are $\binom{N}{n}$ possible ways to choose the distinct $a_k$'s and $n!$ ways to permute them.
We thus obtain
$$Tr(P_n^N)=\frac{1}{n!}\times\binom{N}{n}\times n!$$
So, the value of the trace is simply the binomial coefficient $\binom{N}{n}$.
\begin{equation} Tr(P_n^N)=\binom{N}{n}\end{equation}
Let the $c_i$'s denote the disjoint cycles composing $\sigma$, with their respective lengths $l(c_i)$.\\
We have
\begin{equation}\varepsilon(\sigma)=\prod_{i=1}^{s(\sigma)}\varepsilon(c_i)=\prod_{i=1}^{s(\sigma)}(-1)^{l(c_i)-1}=(-1)^{\sum_{i=1}^{s(\sigma)}(l(c_i)-1)}=(-1)^{n-s(\sigma)}\end{equation}
By gathering Eqs. $(4)$, $(5)$ and $(6)$, we thus get
\begin{equation}
N(N-1)\dots (N-n+1)=\sum_{\sigma\in S_n}(-1)^{n-s(\sigma)}\,N^{s(\sigma)}
\end{equation}
\begin{Definition} The unsigned Stirling number of the first kind $\left[ n\atop s\right]$ is the unsigned coefficient of $x^s$ in the falling factorial
$$(x)_n=x(x-1)\dots(x-n+1)$$
That is, $$\left[n\atop s\right]=\sigma_{n-s}\,(0,\,1,\,\dots,\,n-1),$$ where $\sigma_{n-s}$ denotes the $(n-s)$-th elementary symmetric polynomial.
\end{Definition}
\noindent So, we also have
\begin{equation}
N(N-1)\dots (N-n+1)=\sum_{s=1}^n \left[n\atop s\right]\,(-1)^{n-s} N^s
\end{equation}
With Eqs. $(7)$ and $(8)$, we recover another definition for the unsigned Stirling number of the first kind.
\begin{Definition} The unsigned Stirling number of the first kind $\left[n\atop s\right]$ counts the number of permutations of $n$ elements with $s$ disjoint cycles.
\end{Definition}
\noindent For the interest of the reader, we recall here the more classical way of seeing the equivalence between Definition $3$ and Definition $4$. It can easily be seen from the following recursion formula, which holds with both definitions.
\begin{equation}
\left[n+1\atop s+1\right]=\left[n\atop s\right]+n\,\left[n\atop s+1\right]
\end{equation}
Moreover, both definitions benefit from the same initial conditions, which solves the recursion uniquely. More details can be found on Wikipedia. We slightly continue the digression by giving an analog of the binomial formula with the unsigned Stirling numbers of the first kind. \\
Recall the binomial formula
\begin{equation}
(N-1)^n=\sum_{s=0}^n\, \binom{n}{s}\, (-1)^{n-s}\, N^s
\end{equation}
Our similar formula reads
\newtheorem{Proposition}{Proposition}
\begin{Proposition}
\begin{equation}\binom{N}{n}=\sum_{s=0}^n\, F_s^n\, (-1)^{n-s}\, N^s,\end{equation} with
\begin{equation}F_s^n=\sum_{\begin{array}{l}\text{Ferrers diagrams}\\ \text{partitioning $n$}\\ \text{with $s$ rows}\end{array}}\frac{1}{\prod_{i=1}^s l(R_i)}=\;\sum_{\begin{array}{l}\la_1+\dots +\la_s=n\\\la_1\geq \la_2\geq\dots\geq\la_s\end{array}}\,\frac{1}{\la_1\dots\la_s}\end{equation}
and $l(R_i)$ the length of Row $i$. \\\\The $F_s^n$'s satisfy to the recursion formula
\begin{equation}(n+1)\,F^{n+1}_{s+1}=F_s^n+n\,F_{s+1}^n\end{equation}
\end{Proposition}
\noindent\textbf{Proof.} First, the two different ways of computing the traces resulting in Eq. $(4)$ and Eq. $(5)$ respectively or directly Eq. $(8)$ provides us with the equality
\begin{equation}
\binom{N}{n}=\frac{1}{n!}\,\sum_{s=1}^n\,(-1)^{n-s}\,N^s\,\left[n\atop s\right]\end{equation}
where $\left[n\atop s\right]$ is the number of permutations of $Sym(n)$ that decompose into a product of $s$ cycles with disjoint supports. Since such a decomposition is unique up to commuting its factors, we may always order the cycles composing the product in decreasing order of lengths. Thus, there are obviously $n!\,F_s^n$ such permutations, where the denominator of $F^n_s$ comes from the fact that there are $l$ different ways of writing a cycle of length $l$. It yields Eq. $(11)$. \\
The recursion formula is simply obtained by using Eq. $(9)$ and
$\left[n\atop s\right]=n!\,F_s^n$\hfill$\square$\\

We now use the discussion with the trace to see how the nice specialization $A=-1$ can be used heuristically for computing a closed network of projectors. When $A=-1$, we have $d=-2$, so that by Theorem $1$ and its proof, we have
\begin{center}
\epsfig{file=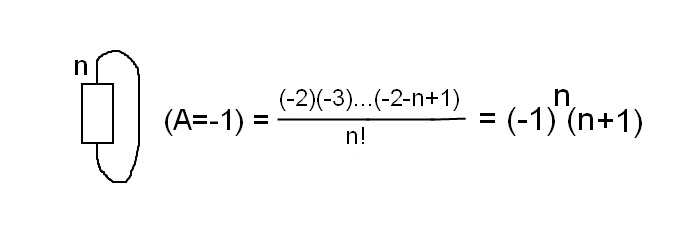, height=4cm}
\end{center}
Note that the first equality holds because once we know the equivalence between Definition $3$ and Definition $4$, 
the equality
$$\sum_{\sigma\in S_n}\,\varepsilon(\sigma)\,N^{s(\sigma)}=N(N-1)\dots (N-n+1)$$ holds for any $N$, and not only for $N$ a positive integer, hence in particular holds for $N=-2$. \\
On the other hand, we know from $\S\,1.2$, Eq. $(3)$ that
\begin{center}
\epsfig{file=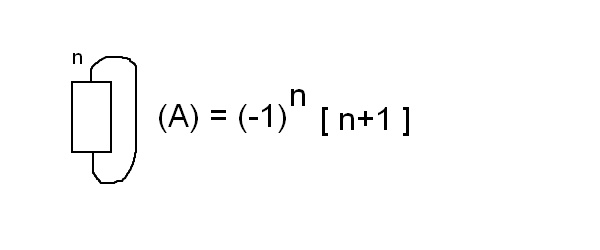, height=4cm}
\end{center}
The idea developed in \cite{KL} is to compute the closed network of projectors at $A=-1$ and then replace all the integers by the quantum integers. \\
We now develop a second heuristics, still arising from the preliminaries.
From now on, say $\lbrace 1,\,\dots,\, N\rbrace$ is a coloring set.
Observe that
\begin{Lemma}
$$[P_n^N]_{b_1\dots b_n}^{a_1\dots a_n}=\begin{cases} \frac{1}{n!}\,\varepsilon(\sigma)&\begin{array}{l}\text{if the $a_i$'s are $n$ distinct colors chosen among $N$ colors}\\\text{and the $b_i$'s are permutations of the $a_i$'s}\\\text{such that $b_{\sigma(i)}=a_i$}\end{array}\\
0&\;\;\text{otherwise}\end{cases}$$
\end{Lemma}
\noindent\textbf{Proof.} Use Lemma $2$ and the definition of $P_n^N$. \hfill $\square$ \\\\
This observation and the discussion from the preliminaries tell us that dealing with tangles and counting the number of loops is the same as dealing with matrix coefficients and summing over all admissible colorings. An admissible coloring is such that the strands are labeled with colors and each loop carries a different color.\\
We use these two heuristics in the next section.

\subsection{Evaluation of the tetrahedron}
We want to evaluate $T\left[\begin{array}{ccc}a&b&e\\c&d&f\end{array}\right]$. \\\\
\begin{center}
\epsfig{file=tetrahedron.jpg, height=4cm}
\end{center}

\noindent For each of the four nodes, there are three values for the internal edges.\\

$$\begin{array}{ccc}\text{Set}&&\left|\begin{array}{ccc}
i=\frac{a+b-f}{2}&j=\frac{b+f-a}{2}&k=\frac{a+f-b}{2}\\
&&\\
l=\frac{a+d-e}{2}&m=\frac{a+e-d}{2}&n=\frac{e+d-a}{2}\\
&&\\
o=\frac{d+c-f}{2}&p=\frac{f+c-d}{2}&q=\frac{f+d-c}{2}\\
&&\\
r=\frac{b+c-e}{2}&s=\frac{b+e-c}{2}&t=\frac{e+c-b}{2}
\end{array}\right.\end{array}$$

\noindent After replacing each node with three projectors, here is the diagram of projectors that we get.
\begin{center}
\epsfig{file=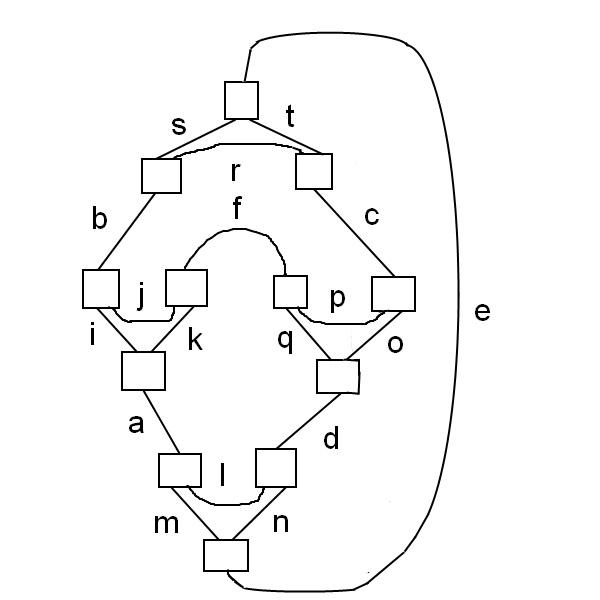, height=9cm}
\end{center}
This diagram simplifies to
\begin{center}
\epsfig{file=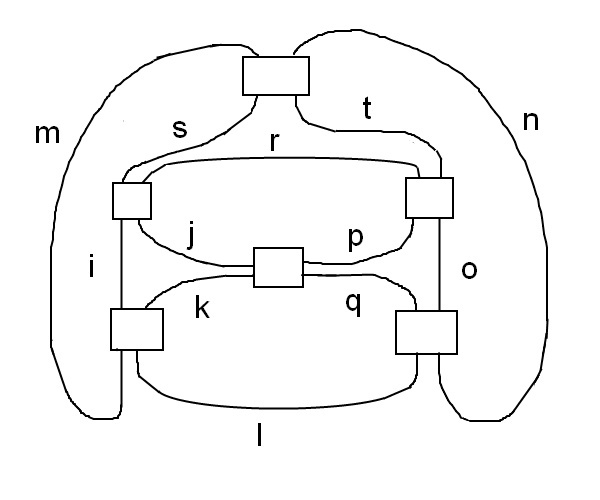, height=8cm}
\end{center}
It will be useful to keep track of the following equalities
$$\left\lbrace\begin{array}{ccccc}
i+j&=&r+s&=&b\\
i+k&=&l+m&=&a\\
j+k&=&p+q&=&f\\
p+o&=&t+r&=&c\\
o+q&=&l+n&=&d\\
s+t&=&m+n&=&e
\end{array}\right.$$
and to introduce notations for the sum of the values of the internal edges.
$$\begin{array}{cc}\text{Define}&\left|\begin{array}{ccccc} T_1&=&i+j+k&=&\frac{a+b+f}{2}\\
&&&&\\
T_2&=&l+m+n&=&\frac{a+d+e}{2}\\
&&&&\\
T_3&=&o+p+q&=&\frac{d+c+f}{2}\\
&&&&\\
T_4&=&r+s+t&=&\frac{b+c+e}{2}
\end{array}\right.\end{array}$$
Other values below will be relevant.
$$\left|\begin{array}{ccccccc}
S_1&=&T_1+T_3-f&=&T_2+T_4-e&=&\frac{a+b+c+d}{2}\\
&&&&&&\\
S_2&=&T_1+T_4-b&=&T_2+T_3-d&=&\frac{a+c+e+f}{2}\\
&&&&&\\
S_3&=&T_3+T_4-c&=&T_1+T_2-a&=&\frac{b+d+e+f}{2}
\end{array}\right.$$
Without loss of generality, we will assume that $r$ is the smallest value of all the values of the internal edges. \\\\
We must compute
\begin{large}
$$\begin{array}{l}
[P^N_b]^{\al_1\dots\al_s\,\al_{s+1}\dots\al_{s+r}}_{\be_1\dots\be_i\,\be_{i+1}\dots\be_{i+j}}\times[P^N_a]^{\be_1
\dots\be_i\,\ga_1\dots\ga_k}_{\de_1\dots\de_m\,\et_1\dots\et_l}
\times [P^N_f]^{\ga_k\dots\ga_1\,\be_{i+1}\dots\be_{i+j}}_{\xi_1\dots\xi_q\,\ep_1\dots\ep_p}\\\\\times [P^N_c]^{\nu_1\dots\nu_o\,\ep_1\dots\ep_p}_{\tau_1\dots\tau_t\,\al_{s+1}\dots\al_{s+r}}
\times[P^N_d]^{\xi_1\dots\xi_q\,\nu_o\dots\nu_1}_{\et_l\dots\et_1\,\ka_1\dots\ka_n}
\times[P^N_e]^{\de_1\dots\de_m\,\ka_1\dots\ka_n}_{\al_1\dots\al_s\,\tau_t\dots\tau_1}\end{array}
$$\end{large}
\noindent where the greek letters at the top denote distinct colors within a same row and the greek letters at the bottom are permutations of those at the top, and we must sum over all the possible colorings. A quick glance at the expression above indicates that there are at most $r+s+k+o$ colors involved.\\
Start from a color $\al_{s+w}$ with $1\leq w\leq r$. It is straightforward to see that there are only two possibilities.
\begin{center}
\epsfig{file=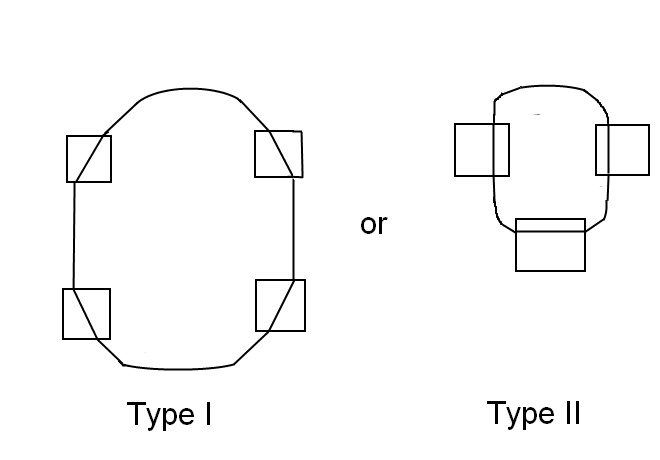, height=7cm}
\end{center}
\newcounter{counter}
\setcounter{counter}{1}
Choose $s+r+k$ distinct colors. The number of colors imposed by the $\al_{s+w}$'s on $\nu_1$, $\dots$, $\nu_o$ equals $|\mcl_{\Roman{counter}}|$, the number of loops of type \Roman{counter}. So, the number of colors to choose for $\nu_1$, $\dots$, $\nu_o$ equals $o-|\mcl_{\Roman{counter}}|$. Thus, the
total number of colors to choose equals $r+s+k+o-\lone$. We name this quantity $Q$. Once these colors are chosen, there are no new colors involved. We next discuss the
minimal and maximal values which $Q$ may take. \\
Since $r\leq p$ and $r\leq j$, we may have all the loops of type \setcounter{counter}{2}\Roman{counter}, in which case \setcounter{counter}{1}$\lone=0$. Then, the number of colors to choose is
$$i+j+k+o=S_1$$
Moreover, since $r\leq p$ and $r\leq j$, we have $r-p\leq p-r$ and $r-j\leq j-r$. But also,
$p-r=t-o$ and $j-r=s-i$. Hence we derive
$$\left\lbrace\begin{array}{l}r-p\leq t-o\\r-j\leq s-i\end{array}\right.$$
Equivalently we have $$\left\lbrace\begin{array}{l}r-t\leq p-o\\r-s\leq j-i\end{array}\right. ,$$ and so $$\left\lbrace\begin{array}{l}b-e\leq f-d\\c-e\leq f-a\end{array}\right.$$ We deduce the following set of inequalities.
$$\left\lbrace\begin{array}{ccc}
a+b+c+d&\leq&a+c+e+f\\
a+b+c+d&\leq&b+d+e+f
\end{array}\right.$$

\noindent We conclude that the maximal number of colors to choose is $S_1=Min_{1\leq w\leq 3}S_w$.
\setcounter{counter}{2}

\noindent Conversely, the minimal number of colors to choose is obtained when \setcounter{counter}{1}$\lone=r$. For that, the following inequalities must hold.
$$\begin{array}{ccccccc}
i\geq r&\Lra& a-f\geq c-e&\Lra&a+d+e\geq d+c+f&\Lra&T_2\geq T_3\\
&&&&&&\\
l\geq r&\Lra&a+d\geq b+c&\Lra&a+d+e\geq b+c+e&\Lra& T_2\geq T_4\\
&&&&&&\\
o\geq r&\Lra&d-f\geq b-e&\Lra& d+e+a\geq b+f+a&\Lra&T_2\geq T_1
\end{array}$$
So, we must have $$T_2=Max_{1\leq x\leq 4}\,T_x$$
This condition holds since $r$ is the smallest label for the internal edges.\\ Thus, the minimal number of colors to choose is
$$r+s+k+o-r=s+k+o=\frac{a+d+e}{2}=T_2=Max_{1\leq x\leq 4}\,T_x$$
Hence the proposition.
\begin{Proposition}
The number $Q$ of colors to choose can vary between $Max_{1\leq x\leq 4}\,T_x$ and $Min_{1\leq w\leq 3}\,S_w$ and is $i+j+k+o-\lone=Q$.
\end{Proposition}
Next, we must choose $\lone$ colors among $r$ colors and there are
$$\frac{r!}{(S_1-Q)!(Q-T_2)!}$$ ways to do so. Then, we must choose $(i-\lone)$ colors among $s$ and we can do it in $$\frac{s!}{(Q-T_3)!(S_3-Q)!}$$ different ways. \\
The projector $P^N_b$ can permute the colors in $$i!\,j!$$ different ways. \\
Further, notice the strands labeled by $l$ never carry the colors from the strands labeled by $s$. Hence, we must choose $(l-\lone)$ colors among $k$ colors. This can be done in
$$\frac{k!}{(Q-T_4)!(S_2-Q)!}$$ different ways. Furthermore, the number of ways the projector $P^N_a$ can permute the colors is $$l!\,m!$$
We must still count the different possible endpoints for the strands belonging to loops of type  $\Roman{counter}$ after they have been traversing the projector $P^N_d$. This number is $\binom{o}{\lone}$ and there are $\lone!$ different ways of organizing the colors. The product of these two quantities is
$$\frac{o!}{(o-\lone)!}=\frac{o!}{(Q-T_1)!}$$
It remains to choose how to arrange the colors on the $p$ strands, $q$ strands, $n$ strands and $t$ strands respectively and there are $p!\,q!\,n!\,t!$ possibilities. \\
Recall that the signature of a permutation $\sigma$ is $(-1)^{n_{\sigma}}$, where $n_{\sigma}$ denotes the number of inversions of the permutation $\sigma$. By the way the colored loops are organized, we see that the total number (all the projectors being involved) of color inversions must be even. So, the final result we get is
$$\frac{\underset{\begin{array}{l}1\leq w\leq 3\\1\leq x\leq 4\end{array}}{\prod}\,(S_w-T_x)!}{a!\,b!\,c!\,d!\,e!\,f!}\;\sum_{Q=\underset{1\leq x\leq 4}{Max}\, (T_x)}^{\underset{1\leq w\leq 3}{Min}\,(S_w)}\,\frac{N(N-1)\dots(N-Q+1)}{\underset{1\leq x\leq 4}{\prod}\,(Q-T_x)!\underset{1\leq w\leq 3}{\prod}\,(S_w-Q)!},$$
where the nominator inside the sum is the number of ways to choose the $Q$ colors and to permute them (this sets the types and positions of the $Q$ colors), the nominator outside the sum is the product of the factorials of the values of the internal edges and the denominator outside the sum corresponds to the factorials appearing in our Lemma $3$ (one per projector). By our second heuristics, this quantity is the evaluation of the network at $d=N$ when $A=-1$. When doing $d=-2$, we get
$$\frac{\underset{\begin{array}{l}1\leq w\leq 3\\1\leq x\leq 4\end{array}}{\prod}\,(S_w-T_x)!}{a!\,b!\,c!\,d!\,e!\,f!}\;\sum_{Q=\underset{1\leq x\leq 4}{Max}\, (T_x)}^{\underset{1\leq w\leq 3}{Min}\,(S_w)}\,\frac{(-1)^Q(Q+1)!}{\underset{1\leq x\leq 4}{\prod}\,(Q-T_x)!\underset{1\leq w\leq 3}{\prod}\,(S_w-Q)!},$$
In order to obtain the evaluation of the network for a generic $A$, by the first heuristics, it suffices to replace all the integers by the quantum integers taken at $q=A^2$. We obtain the formula
$$T\left[\begin{array}{ccc}a&b&e\\c&d&f\end{array}\right]=\frac{\underset{\begin{array}{l}1\leq w\leq 3\\1\leq x\leq 4\end{array}}{\prod}\,[S_w-T_x]!}{[a]!\,[b]!\,[c]!\,[d]!\,[e]!\,[f]!}\;\sum_{Q=\underset{1\leq x\leq 4}{Max} \,(T_x)}^{\underset{1\leq w\leq 3}{Min}\,(S_w)}\,\frac{(-1)^Q[Q+1]!}{\underset{1\leq x\leq 4}{\prod}\,[Q-T_x]!\underset{1\leq w\leq 3}{\prod}\,[S_w-Q]!}$$
\subsection{Discussion}
This general formula presented here in a heuristic approach different from \cite{KL} was proven rigourously in \cite{MV} by induction using the Jones-Wenzl recursion relation.
Induction is safe and systematic, but has nothing of a light and pleasant method to get to the result. Thus, we ask whether our heuristic approach could be generalized to evaluating any plane trivalent graph like on the following figure. We answer positively below. We raise however the question whether the heuristics which include the quantization at the final stage will give the correct result. If the answer is yes, we would then have a rigorous method for evaluating generically any plane trivalent graph of Jones-Wenzl projectors, a method which does not use induction. Of course any plane trivalent graph can be evaluated without induction by using a series of moves with the $6j$-symbols (see the recoupling theorem for a generic $A$ in Chapter $7$ of \cite{KL}) and theta and tetrahedron evaluations, but this way seems much harder to program. The latter method is used extensively in \cite{BL} where plane trivalent graphs of small sizes are involved. 
\begin{center}
\epsfig{file=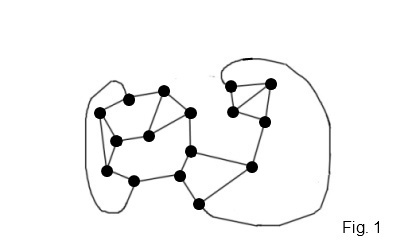, height=6cm}
\end{center}
We now discuss the generalization of our work on the tetrahedron which could apply to any graph like the graph of Fig. $1$ randomly drawn above. First, there is an algorithm for replacing any trivalent plane graph with a graph of Jones-Wenzl projectors. A general method is explained in \cite{KL} at the beginning of $\S\,8.4$. We present here a different way. Notice by the core property of a projector, each projector of the resulting graph will have exactly four outgoing edges which are the so-called internal edges of the original network. Call this resulting graph reduced projector graph. Moreover, each projector of the reduced projector graph will hence be associated with a unique external edge which is an edge from the initial network. Thus, the number $V^{'}$ of vertices of the reduced projector graph equals the number of edges of the trivalent graph, that is
$$V^{'}=\frac{3V}{2},$$ where $V$ denotes the number of vertices of the trivalent graph. It then suffices to form the $4$-valent graph with $\frac{3V}{2}$ vertices and place the internal edges.
Once we have the graph in terms of projectors, here is how we proceed to the evaluation. We can always write the graph as a product of matrix coefficients. A simple reading of this product provides an upper bound for the number of colors to choose. Since each projector takes in entry distinct colors and permutes them, we can determine all the admissible cycles by the only constraint that a color cannot appear twice at an extremity of a projector. For a network of big size in terms of the number of projectors that it contains, a computer program is needed and well suited to determine all the admissible cycles. Next, it is easy to tell the exact number of colors involved. Indeed, when reading the product, we may encounter an edge which will potentially carry new colors. This edge belongs to a certain number of admissible cycles. Then the number of new colors carried by this edge is its label minus the total number of loops of known colors carried by the cycles of the different types which this edge belong to. Without loss of generality, the first matrix coefficient in our product can come from a projector of the network such that one of the labels of its four outgoing edges is a minimal label, say $l$, among all the labels of the graph. The row and column of the first matrix coefficient are then chosen accordingly so that the row contains the colors from these $l$ strands.
From there, the generalization attempt requires far more effort and is the object of ongoing work.
\section{Quantization of the factorial}
From Definition $4$, we have $$n!=\sum_{s=1}^n\left[n\atop s\right]$$
Our goal in this part is to find a similar expression for the quantum factorial
$$[n]!=[n][n-1]\dots [2][1]$$ using the unsigned Stirling numbers of the first kind. It is not easy to quantize a sum in general. Indeed, for non-negative integers $a$ and $b$, we have
$$[a+b]_q=q^b[a]_q+q^{-a}\,[b]_q$$
An idea is to start again from the equality
\begin{equation}N(N-1)\dots(N-n+1)=\sum_{s=1}^n \,(-1)^{n-s}\,N^s\left[n\atop s\right]\end{equation} and to substitute $N=-2$.
Here is the surprising formula that we get.
\begin{equation}(n+1)!=\sum_{s=1}^n\,\left[n\atop s\right]\,2^s\end{equation}
The following proposition gives the quantized version of this equality.
\begin{Proposition}
In the formula below, $\left[n\atop s\right]$ denotes the unsigned Stirling number of the first kind and all the other brackets denote quantum integers. Also, $I\left(\frac{n+1}{2}\right)$ denotes the integer part of $\frac{n+1}{2}$. By definition, $\varepsilon_P$ is $1$ if Proposition $P$ is true and $-1$ otherwise.
\begin{equation}[n+1]!=\sum_{s=1}^n\,\left[n\atop s\right][2]^s+
\la_0^{(n+1)}+\sum_{l=1}^{\binom{n+1}{2}}\la_l^{(n+1)}\,\frac{[2l]}{[l]}\end{equation}
The integers $\la_i^{(n+1)}$'s with $0\leq i\leq \binom{n+1}{2}$ satisfy to
$$\left\lbrace\begin{array}{ccc}
2\,\underset{\begin{array}{l}1\leq l\leq \binom{n+1}{2}\\\text{$l$ is even}\\\text{$\frac{l}{2}$ is odd}\end{array}}{\sum}\,\la_l^{(n+1)}&=&\frac{(n+1)!}{4}\,\varepsilon_{\lbrace I\left(\frac{n+1}{2}\right)\;\text{is even}\rbrace}\\&&\\\la_0^{(n+1)}+2\,\underset{\begin{array}{l}1\leq l\leq \binom{n+1}{2}\\\text{$l$ is even}\\\text{$\frac{l}{2}$ is even}\end{array}}{\sum}\,\la_l^{(n+1)}&=& \frac{(n+1)!}{4}\,\varepsilon_{\lbrace I\left(\frac{n+1}{2}\right)\;\text{is even}\rbrace}\\&&\\2\,\underset{\begin{array}{l} 1\leq l\leq \binom{n+1}{2}\\\text{$l$ is odd}\end{array}}{\sum}\,\la_l^{(n+1)}&=&\frac{(n+1)!}{2}\,\varepsilon_{\lbrace I\left(\frac{n+1}{2}\right)\;\text{is odd}\rbrace}\end{array}\right.$$

\end{Proposition}
\noindent \textbf{Proof.} By definition of the quantum integers, we have
\begin{multline*}
[n+1]!=(q+q^{-1})(q^2+1+q^{-2})(q^3+q+q^{-1}+q^{-3})(q^4+q^2+1+q^{-2}+q^{-4})\\(q^5+q^3+q+q^{-1}+q^{-3}+
q^{-5})\dots (q^n+q^{n-2}+q^{n-4}+\dots + q^{-n+4}+q^{-n+2}+q^{-n})
\end{multline*}
and $$[2]^s=(q+q^{-1})^s$$ and $$\frac{[2l]}{[l]}=q^l+q^{-l}$$
Hence the general formula of Proposition $3$ where the $\la_i^{(n+1)}$'s, $0\leq i\leq \binom{n+1}{2}$, are some integers to determine. In order to get relations between these coefficients, we let $q$ tend to $1$ in Eq. $(17)$ and we use Eq. $(16)$. It yields
\begin{equation}
\la_0^{(n+1)}+2\,\sum_{l=1}^{\binom{n+1}{2}}\,\la_l^{(n+1)}=0
\end{equation}
We further let $q$ tend to $-1$ in Eq. $(17)$ and thus get
\begin{multline}
\la_0^{(n+1)}+2\,\underset{\begin{array}{l}1\leq l\leq \binom{n+1}{2}\\\text{$l$ is even}\end{array}}{\sum}\,\la_l^{(n+1)}-2\,\underset{\begin{array}{l}1\leq l\leq \binom{n+1}{2}\\\text{$l$ is odd}\end{array}}{\sum}\,\la_l^{(n+1)}\\+(-1)^n\,\sum_{s=1}^n\,\left[n\atop s\right](-1)^{n-s}\,2^s=(-1)^{I\left(\frac{n+1}{2}\right)}\,(n+1)!
\end{multline}
By Eq. $(15)$ with $N=2$, the third sum in Eq. $(19)$ is zero. So, Eq. $(19)$ simplifies to
\begin{equation}\la_0^{(n+1)}+2\,\underset{\begin{array}{l}1\leq l\leq \binom{n+1}{2}\\\text{$l$ is even}\end{array}}{\sum}\,\la_l^{(n+1)}-2\,\underset{\begin{array}{l}1\leq l\leq \binom{n+1}{2}\\\text{$l$ is odd}\end{array}}{\sum}\,\la_l^{(n+1)}=(-1)^{I\left(\frac{n+1}{2}\right)}\,(n+1)!\end{equation}
From Eq. $(18)$ and Eq. $(20)$, we derive
$$\begin{array}{ccc} \underset{\begin{array}{l}1\leq l\leq \binom{n+1}{2}\\\text{$l$ is odd}\end{array}}{\sum}\,\la_l^{(n+1)}&=&\begin{cases} -\frac{(n+1)!}{4}&\text{if $I\left(\frac{n+1}{2}\right)$ is even}\\\frac{(n+1)!}{4}&\text{if $I\left(\frac{n+1}{2}\right)$ is odd}\end{cases}\\
\la_0^{(n+1)}+2\,\underset{\begin{array}{l}1\leq l\leq \binom{n+1}{2}\\\text{$l$ is even}\end{array}}{\sum}\,\la_l^{(n+1)}&=&\begin{cases}\frac{(n+1)!}{2}&\text{if $I\left(\frac{n+1}{2}\right)$ is even}\\
-\frac{(n+1)!}{2}&\text{if $I\left(\frac{n+1}{2}\right)$ is odd}\end{cases}
\end{array}$$
Furthermore, by substituting $q=i$ in Eq. $(17)$, we have
\begin{equation}
\la_0^{(n+1)}+2\,\underset{\begin{array}{l}1\leq l\leq \binom{n+1}{2}\\\text{$l$ is even}\end{array}}{\sum}\,\la_l^{(n+1)}\,(-1)^{\frac{l}{2}}=0
\end{equation}
So, we get \begin{equation}\underset{\begin{array}{l}1\leq l\leq \binom{n+1}{2}\\\text{$l$ is even}\\\text{$\frac{l}{2}$ is odd}\end{array}}{\sum}\,\la_l^{(n+1)}=\frac{(n+1)!}{8}\,\varepsilon_{\lbrace I\left(\frac{n+1}{2}\right)\;\text{is even}\rbrace}\end{equation}
It follows that
\begin{eqnarray*}\la_0^{(n+1)}+2\,\underset{\begin{array}{l}1\leq l\leq \binom{l+1}{2}\\\text{$l$ is even}\\
\text{$\frac{l}{2}$ is even}\end{array}}{\sum}\la_l^{(n+1)}&=&\frac{(n+1)!}{4}\,\varepsilon_{\lbrace I\left(\frac{n+1}{2}\right)\;\text{is even}\rbrace}\\&=&-\underset{\begin{array}{l}1\leq l\leq \binom{l+1}{2}\\\text{$l$ is odd}\end{array}}{\sum}\,\la_l^{(n+1)}\end{eqnarray*}
We obtain the equalities of Proposition $3$. \hfill $\square$
\begin{Proposition} (Recursive determination of the coefficients $\la_i$'s).
Define $$\mu_i^{(n+1)}=\begin{cases}\la_i^{(n+1)}+\underset{\begin{array}{l}i\leq s\leq n\\s\equiv i\,\text{mod}\,2\end{array}}{\sum}\left[n\atop s\right]\,\binom{s}{\frac{s+i}{2}}&\text{when $0\leq i\leq n$}\\\la_i^{(n+1)}&\text{when $n+1\leq i\leq \binom{n+1}{2}$}\end{cases}$$
When $i\geq 1$, $\mu_i^{(n+1)}$ is the coefficient of $(q^i+q^{-i})$ in the expansion of $[n+1]!$ into powers of $q$ and $\mu_0^{(n+1)}$ is the constant coefficient.\\\\
When $I\left(\frac{n+1}{2}\right)$ is odd (resp even), we have
$$\la_i^{(n+1)}=\begin{cases}-\underset{\begin{array}{l}i\leq s\leq n\\s\equiv i\,\text{mod}\,2\end{array}}{\sum}\left[n\atop s\right]\,\binom{s}{\frac{s+i}{2}}&\text{when $i$ is even (resp odd) and $0\leq i\leq n$}\\0&\text{when $i$ is even (resp odd) and $n+1\leq i\leq \binom{n+1}{2}$}\end{cases}$$
Suppose the coefficients $\la_i^{(n+1)}$'s known. Then, the coefficients $\la_i^{(n+2)}$'s are determined recursively by
\begin{multline*}\mu_i^{(n+2)}=\mu_{|n+1-i|}^{(n+1)}+\,\mu_{|n-1-i|}^{(n+1)}+\dots+\mu_{|\delta_{\lbrace\text{$n$
is even}\rbrace}-i|}^{(n+1)}\\+\mu_{|n+1+i|}^{(n+1)}+\,\mu_{|n-1+i|}^{(n+1)}+\dots+\mu_{|\delta_{\lbrace\text{$n$
is even}\rbrace}+i|}^{(n+1)}\,\delta_{\lbrace\text{$n$ is even}\rbrace}\end{multline*}
We ask whether it is possible to find explicit formulas for the coefficients $\la_i$'s.
\end{Proposition}
\noindent \textbf{Proof.} For the second point, it suffices to notice that when $I\left(\frac{n+1}{2}\right)$ is odd (resp even), the exponents of the powers of $q$ in the expansion of $[n+1]!$ are all odd (resp even). The rest is clear. \hfill $\square$\\\\
We now state an a priori surprising fact.
\begin{Proposition}
Write $$[n+1]!=(q+q^{-1})(q^2+1+q^{-2})\dots (q^n+q^{n-2}+\dots+q^{-n+2}+q^{-n})=\sum_{k=-\binom{n+1}{2}}^{\binom{n+1}{2}} c_k\,q^k$$ and
$$\sum_{s=1}^n\left[ n\atop s\right]\,[2]^s=\left[n\atop 1\right](q+q^{-1})+\left[n\atop 2\right](q+q^{-1})^2+\dots + \left[n\atop n\right](q+q^{-1})^n=\sum_{k=-n}^n\,d_k\,q^k$$
We have $$\sum_{k=-\binom{n+1}{2}}^{\binom{n+1}{2}}\,c_k=\sum_{k=-n}^n\,d_k=(n+1)!$$
In summary, $$(n+1)!=\sum_{s=1}^n\left[n\atop s\right]2^s$$
$$[n+1]!\neq \sum_{s=1}^n\left[n\atop s\right][2]^s,$$
but the sum of the coefficients in the respective Laurent series are equal to $(n+1)!$ 
\end{Proposition}
\noindent\textbf{Proof.} Follows from Eq. $(17)$ and Eq. $(18)$ or simply notice that these sums are obtained by letting $q$ tend to $1$, hence using Eq. $(16)$ are in both cases $(n+1)!$ \\
\noindent We illustrate the result on an example. We have
\begin{multline*}[5]!=(q^{10}+q^{-10})+4\,(q^8+q^{-8})+9\,(q^6+q^{-6})\\+15(q^4+q^{-4})+20(q^2+q^{-2})+22
\end{multline*}
\begin{multline*}\sum_{s=1}^4\left[4\atop s\right](q+q^{-1})^s=
(q^4+q^{-4})+6(q^3+q^{-3})+15(q^2+q^{-2})+24(q+q^{-1})+28\end{multline*}
and as a matter of fact,
$$28+48+30+12+2=22+40+30+18+8+2$$
Observe in the sum with the Stirling numbers, the sum of the coefficients of the even powers of $q$ equals the sum of the coefficients of the odd powers of $q$.
$$28+30+2=48+12$$This fact is not random and actually holds for all $n$. Indeed, we have the following result.
\begin{Proposition}
In $\sum_{s=1}^n\left[n\atop s\right][2]^s$, the sum $\mathcal{S}_e$ of the coefficients of the even powers of $q$ equals the sum $\mathcal{S}_o$ of the coefficients of the odd powers of $q$. Moreover, this sum is
$$\sum_{k=1}^{I(\frac{n}{2})}\,\left[ n\atop 2k\right]\,2^{2k}=\sum_{k=1}^{I\left(\frac{n-1}{2}\right)+1}\,\left[n\atop 2k-1\right]\,2^{2k-1}=\frac{(n+1)!}{2}$$
\end{Proposition}
\noindent \textbf{Proof.}
By specializing $q=-1$, we get $$\mathcal{S}_e-\mathcal{S}_o=\sum_{s=1}^n\,(-1)^s\,2^s\left[n\atop s\right]$$
We conclude that $\mathcal{S}_e=\mathcal{S}_o$ by noticing that the right hand side is zero after specializing $N=2$ in Eq. $(15)$.
By specializing $q=1$, we also have
$$\mathcal{S}_e+\mathcal{S}_o=\sum_{s=1}^n\,2^s\left[n\atop s\right]$$
We thus derive \begin{equation}\mathcal{S}_o=\sum_{\begin{array}{l} 1\leq s\leq n\\\text{$s$ is odd}\end{array}}\left[n\atop s\right]2^s=\sum_{\begin{array}{l} 1\leq s\leq n\\\text{$s$ is even}\end{array}}\left[n\atop s\right]2^s=\mathcal{S}_e=\frac{(n+1)!}{2}\end{equation}
Note that the middle equality in $(23)$ between the weighted sums of unsigned Stirling numbers of the first kind also holds when $2$ is being replaced with any integer $k$ such that $1\leq k\leq n-1$. This fact is easily seen by simply substituting $N=k$ in Eq. $(15)$.
When $k=1$, the equality says that the number of permutations on $n$ letters with an odd number of cycles equals the number of permutations on $n$ letters with an even number of cycles. We ask whether there is any combinatorial interpretation of these equalities for the other values of $k$. \\\\

\noindent We would like to finish this paper by underlining the beauty of the algebraic and quantum combinatorics which allow to get rid of induction, like for instance when we provide two equivalent definitions for the unsigned Stirling numbers of the first kind or when we obtain a formula for a tetrahedral network of Jones-Wenzl projectors.

\noindent \\\textbf{Acknowledgements.} The author thanks Zhenghan Wang for introducing her to the subject with much kindness and time generosity and for his comments on this manuscript. She thanks Michael Freedman for helpful discussions. She is grateful to the mathematics department of the University of California at Santa Barbara, Microsoft Research Station Q and Harish Chandra Research Institute where this work was produced for very kind hospitality.

\end{document}